\documentclass{article}%
\usepackage{amsmath}
\usepackage{amsfonts}
\usepackage{amsthm}
\usepackage{lineno,hyperref}
\usepackage{amssymb}
\usepackage{graphicx}
\usepackage{epstopdf}
\usepackage{caption}
\usepackage{subcaption}%
\usepackage{natbib}
\usepackage{authblk}
\newtheorem{theorem}{Theorem}

\begin{document}
\vfill
{\centering%
\noindent\begin{tabular}{|l|}
\hline
\textit{Sanz, L., 2019. Conditions for growth and extinction in matrix models}\\
\textit{ with environmental stochasticity, Ecological Modelling, 411, p.108797.}\\
\textit{https://doi.org/10.1016/j.ecolmodel.2019.108797}\\
\hline
\end{tabular}
}%
\vfill

\title{Conditions for growth and extinction in matrix models with environmental stochasticity}
\author[1]{Luis Sanz\thanks{luis.sanz@upm.es}}
\affil[1]{Depto. Matem\'aticas, E.T.S.I Industriales, Technical University of Madrid, Madrid, Spain.}
\date{}

\maketitle

\begin{abstract}
In this kind of model, the main characteristic that determines population viability in the long term is the \textit{stochastic growth rate} (SGR) denoted $\lambda_S$. When $\lambda_S$ is larger
than one, the population grows exponentially with probability one and when
it is smaller than one, the population goes extinct with probability
one. However, even in very simple situations it is not possible to calculate
the SGR analytically. The literature offers some approximations for the
case in which environmental variability is low, and there are also some lower
and upper bounds, but there is no study of the practical situations in which they would be tight.
Some new bounds for the SGR are built and the conditions under which each
bound works best are analyzed. These bounds are used to give some necessary and some sufficient conditions for population explosion and extinction that are easy to check in practice. The
general results are applied to several cases, amongst them a population structured as juveniles and adults living in an environment switching randomly between "rich" and "poor".
\end{abstract}
%
{\bf keywords:}matrix models, environmental stochasticity, stochastic growth rate, lower bound, upper bound, extinction.



\section{Introduction}

We consider matrix population models that incorporate environmental
stochasticity, which refers to unpredictable temporal fluctuations in
environmental conditions. A number of different environmental conditions are
considered and their variation is characterized by a sequence of random
variables corresponding to the different time steps of the discrete system
\citep{Tulja90,Caswell01}. Each environment is characterized by its
corresponding population projection matrix (PPM), i.e., the matrix of vital
rates in that environment, e.g., time-specific annual PPMs gained from a
time-series of data (Logofet, 2018; elsewhere in this Issue). It is well known
that, given certain hypotheses on the pattern of temporal variation and on the
set of PPMs, the total population size grows (or decays) exponentially with
probability one, with a rate given by the so called \textit{stochastic growth
rate} $\lambda_{\mathrm{S}}$ (SGR), which is the stochastic analogue of the
dominant eigenvalue for deterministic systems.

The SGR is the main feature of this kind of models due to its ecological
implications, and it is widely used by ecologists as a measure of population
fitness
\citep{williams2015life,compagnoni2016effect,mcdonald2016transients,pmid:24393387}.
A SGR larger than one implies population explosion with probability one
whereas a value lower than one implies that the population tends to extinction
with probability one. It must be stressed that we are adopting Cohen's
definition for the SGR \citep{cohen1979long} because it allows one to directly
compare the stochastic rate of growth with its deterministic counterpart,
although some authors \citep{Tulja90,Caswell01} define the SGR as the
logarithmic rate of growth, i.e., $\log\lambda_{\mathrm{S}}$.

The SGR is easy to obtain for scalar models, i.e., in the case of unstructured
populations \citep{alonso2009approximating}. However (letting aside some
settings in which very specific conditions are met
\citep{Tulja90,alonso2009approximating}), it can not be derived analytically
for structured populations. Indeed, even in relatively simple cases like
stochastic Leslie models with two age classes, there are no closed-form
expressions for the SGR.

As a result, in order to obtain the SGR, the main tools at the ecologist's
disposal are computer simulations\citep{williams2015life,pmid:24393387},
certain perturbation techniques like Tuljapurkar's approximation, $\lambda
_{T},$ to the SGR (see \citep{Tulja90} for the derivation and
\citep{davison2013contributions} for an application) and, in certain systems
with time scales, approximate aggregation techniques that allow one to
transform the original system to a scalar one \citep{sanz2007approximate,alonso2009approximating,sanz2010approximate}.

Another strategy that offers valuable information is the derivation of upper
and lower bounds for the SGR. In addition to the quantitative values of these
bounds, ecologists can make use of them to extract qualitative information
about the model in the following way: assume that one knows that%
\begin{equation}
c\leq\log\lambda_{\mathrm{S}}\leq C, \label{qaz}%
\end{equation}
(throughout this work, and for the sake of simplicity in the expressions
involved, bounds for $\log\lambda_{\mathrm{S}},$ rather than for
$\lambda_{\mathrm{S}}$ itself, will be considered. Obviously bounds for
$\lambda_{\mathrm{S}}$ follow in the form $\exp(c)\leq\lambda_{\mathrm{S}}%
\leq\exp(C)$). Then, the population grows to infinity with probability one
when $c>0$, whereas whenever $C<0$ the population goes extinct with
probability one. Therefore, the bounds can be used to obtain sufficient
(although not necessary) conditions for population growth or extinction.

In the mathematical literature, there are some bounds for the SGR\ that can be
quite precise, but they must be obtained through complex algoritmic
procedures, so that there are no closed-form expressions for them (see
\citep{protasov2013lower} and the references therein). In this work, however,
we are interested in bounds that can be easily obtained and that, at least in
two-stage models, can be obtained analytically. The literature on mathematical
ecology presents a few bounds with those properties
\citep{cohen1979long,Caswell01}. One of them is $\mu$, the rate of growth of
the mean population size, which is an upper bound. Several works have compared
$\lambda_{\mathrm{S}}$ and $\mu$ (see \citep{Tulja90} and the references
therein). However, there is hardly any study on the conditions under which the
other bounds work well or a comparison between the different bounds. An
exception was given by \citep{tuljapurkar1980population}, who discussed a case
in which one of the bounds is very loose and therefore of little practical
interest. Neither are there studies that use the bounds to obtain analytical
conditions for growth or extinction as suggested above.

The first purpose of this work is to find new upper and lower bounds for the
SGR. Different approaches are followed and, as a consequence, four different
upper and lower bounds will be obtained. One of the approaches is only valid
in the context of stochastic Leslie models with two age classes, but the rest
of the bounds are valid for general stage-classified models. Moreover, in the
case of models with two stages, closed-form expressions for the bounds are provided.

Our second aim is to provide some properties of both the new and known bounds
and, moreover, to give some ideas about the situations in which each one
performs best. In addition, a comparison between the different bounds is
carried out.

Finally, in order to illustrate their usefulness in a practical application,
the bounds are used to obtain conditions that guarantee non-extinction in a
Leslie model with two age classes and two possible environments.

The structure of the paper is as follows: Section \ref{sec:matrixmodels}
introduces the basic setting of matrix models with environmental stochasticity
and the concept of SGR. The general setting is particularized in
\ref{sec:leslie} to a stochastic Leslie model with two age classes that will
be used throughout the paper to illustrate the different results we provide.
Section \ref{sec:knownbounds} presents the different bounds already available
in the literature and some of their properties. In Section \ref{sec:pert}, we
make use of a perturbation approach to obtain a new result that is valid for
any kind of stochastic matrix model. It allows one to compare the SGR of two
stochastic systems that share the same pattern of environmental variability
but have different PPMs. Then we make use of this general result in order to
obtain three new bounds for the SGR.

Section \ref{sec:boundsagestructure} concentrates in the specific case of
stochastic Leslie models with two age classes, and provides bounds that are
built by making use of some properties of age structure.

An analysis and comparison of the different bounds is carried out in Section
\ref{sec:comparison} in the context of the aforementioned stochastic Leslie
model with young and adults and two environments. We also relate the bounds
with Tuljapurkar's second-order approximation $\lambda_{T}$ to the SGR \citep{tuljapurkar1982population}.

Section \ref{sec:delta}, towards which ecologists interested mainly in
applications can focus their attention, illustrates one of the possible uses
of the bounds to obtain analytical conditions for growth or extinction of
populations in a setting in which a population is subjected to a rich and a
poor environment. A discussion of results can be found in Section
\ref{sec:discussion}. The manuscript is completed with an Appendix that
contains different proofs together with an important auxiliary result on
stochastic matrix population models.

\section{Matrix models with environmental
stochasticity\label{sec:matrixmodels}}

We assume that the population lives in an ambient in which there are different
environmental states, which we suppose indexed by $\eta\in\mathcal{I}$ where
$\mathcal{I}$ is a set of indexes (examples are $\mathcal{I}=\left\{
1,...,r\right\}  $, $\mathcal{I}=\mathbb{Z}_{+}$ or $\mathcal{I}%
=\mathbb{R}_{+}$). The vital rates corresponding to environment $\eta
\in\mathcal{I}$ are given by the nonnegative PPM $\boldsymbol{A}_{\eta
}=\left[  A_{\eta}^{ij}\right]  _{1\leq i,j\leq n}\in\mathbb{R}_{+}^{n\times
n}$ and $\mathcal{A}=\{\boldsymbol{A}_{\eta},\ \eta\in\mathcal{I}\}$ is the
set of different PPMs.

Note that $\mathcal{I}$ can theoretically be infinite, numerable or
denumerable, yet for the sake of simplicity we express all expectations
considering the numerable case, i.e., in terms of (possible infinite) sums. In
the denumerable case all expressions are still valid but the sums must be
replaced with the corresponding integrals.

Environmental randomness is characterized by a homogeneous Markov chain
$\tau_{t}$, $t=0,1,2,...$ defined on a certain probability space
$(\Omega,\mathcal{F},p)$ \citep{billingsley2008probability} over the state
space $\mathcal{I}$. For each realization $\omega\in\Omega$ of the process,
the population is subjected to environmental conditions $\tau_{t+1}(\omega)$
between times $t$ and $t+1$.

Thus, the model reads
\begin{equation}
\boldsymbol{z}(t+1)=\boldsymbol{A}_{\tau_{t+1}}\boldsymbol{z}(t), \label{e1}%
\end{equation}
where for each $t=0,1,...$, $\boldsymbol{z}(t)=\left(  z^{1}(t),...,z^{n}%
(t)\right)  ^{\text{T}}$ is a vector random variable in $\mathbb{R}^{n}$ which
represents the population vector at time $t$. The total population is
$\left\Vert \boldsymbol{z}(t)\right\Vert _{1}:=\left\vert z^{1}(t)\right\vert
+\cdots+\left\vert z^{n}(t)\right\vert $ (the subscript at the norm notation
will be further omitted). Throughout we assume that $\boldsymbol{z}_{0}$ is a
fixed nonnegative nonzero vector.

Assume that $\tau_{t}$ has a unique stationary probability distribution given
by $\boldsymbol{\pi}=\left[  \pi_{\eta}\right]  _{\eta\in I}$ towards which
the $t$-step transition probabilities converge, i.e., $\pi_{\eta}%
:=\lim_{t\rightarrow\infty}p\left(  \tau_{t}=\eta\mid\tau_{0}=\xi\right)  $
for all $\eta,\xi\in\mathcal{I}$. It is well known (see \ref{app1} for a
precise statement of the result) that if the set $\mathcal{A}$ is ergodic and
an additional technical condition on the speed of convergence of these
transition probabilities is met, then:

\begin{enumerate}
\item The joint process $\left(  \tau_{t},\boldsymbol{z}(t)/\left\Vert
\boldsymbol{z}(t)\right\Vert \right)  $ is a Markov chain that converges to a
certain stationary distribution $G$.

\item There exists the so called \textit{stochastic growth rate}
$\lambda_{\mathrm{S}}$ (SGR) of system (\ref{e1}), defined through
\begin{equation}
\log\lambda_{\mathrm{S}}:=\underset{t\rightarrow\infty}{\lim}\log\left\Vert
\boldsymbol{z}(t)\right\Vert /t, \label{SGR}%
\end{equation}
where the limit holds with probability one. Moreover, $\lambda_{\mathrm{S}}$
is finite and independent of the initial probabilities of the chain states and
of the initial (non-zero) population vector $\boldsymbol{z}_{0}\geq
\boldsymbol{0}$.
\end{enumerate}

It follows from (\ref{SGR}) that if $\log\lambda_{\mathrm{S}}>0,$ then the
population grows exponentially to infinity with probability one, whereas it
becomes extinct with probability one whenever $\log\lambda_{\mathrm{S}}<0$. We
assume the conditions for the existence of the SGR to hold true throughout
this work.

In the next sections several upper and lower bounds for the SGR will be
obtained. Clearly, a desirable property of a SGR bound is that, when there is
no environmental variation (i.e., in the deterministic case), it coincides
with the dominant eigenvalue of the corresponding PPM. Specifically:

\noindent\textit{\textbf{Property P.}} We say that a (lower or upper) bound
$C$ for the SGR verifies property P when in the case $\boldsymbol{A}_{\eta
}=\boldsymbol{A}$ for all$\ \eta\in\mathcal{I}$ it follows that $C=\rho
(\boldsymbol{A}),$ where $\rho$ denotes the spectral radius.

It will be shown that some, but not all, of the bounds studied in this work
verify property P.

\subsection{A particular case:\ stochastic Leslie model with two age
clases\label{sec:leslie}}

As an illustration of the general setting above we turn our attention to the
case of Leslie models with $n=2$ age classes, juveniles and adults. The model
takes the form
\begin{equation}
\left(
\begin{array}
[c]{c}%
\boldsymbol{z}^{1}(t+1)\\
\boldsymbol{z}^{2}(t+1)
\end{array}
\right)  =\boldsymbol{A}_{\tau_{t+1}}\left(
\begin{array}
[c]{c}%
\boldsymbol{z}^{1}(t)\\
\boldsymbol{z}^{2}(t)
\end{array}
\right)  , \label{Leslie}%
\end{equation}
where
\[
\boldsymbol{A}_{\eta}:=\left(
\begin{array}
[c]{cc}%
f_{\eta} & F_{\eta}\\
s_{\eta} & 0
\end{array}
\right)  ,\ \eta\in\mathcal{I}.
\]

This simple model will be used to illustrate the different results obtained.

In model (\ref{Leslie}), a necessary and sufficient condition for the set
$\mathcal{A}$ of vital rates to be ergodic is that there exist positive
constants $\alpha$ and $\beta$ such that for all $\eta\in\mathcal{I}$,
$f_{\eta}$, $F_{\eta}$ and $s_{\eta}$ are bounded from above by $\beta$ and
from below by $\alpha$. In that case, for any non-zero initial condition
$\boldsymbol{z}(0)$ the population vector $\boldsymbol{z}(t)$ is positive for
all $t\geq2$.

The mean matrix with respect to the stationary environmental distribution
$\boldsymbol{\pi}$ is
\begin{equation}
\boldsymbol{\bar{A}}:=\left(
\begin{array}
[c]{cc}%
\bar{f} & \bar{F}\\
\bar{s} & 0
\end{array}
\right)  , \label{LeslieMedia}%
\end{equation}
where
\[
\bar{f}:=\sum_{\eta\in\mathcal{I}}\pi_{\eta}f_{\eta},\ \bar{F}:=\sum_{\eta
\in\mathcal{I}}\pi_{\eta}F_{\eta},\ \bar{s}:=\sum_{\eta\in\mathcal{I}}%
\pi_{\eta}s_{\eta}.
\]

\section{Known bounds for the SGR\label{sec:knownbounds}}

This section shows the bounds already known in the literature. Specifically,
the rate of growth of the mean population size, Cohen's bounds and the bounds
through a max-min approach are considered. Later on, they will be compared
with the new bounds proposed in this work.

\noindent\textit{\textbf{Rate of growth of the mean population size}}. An
important parameter in models of the kind (\ref{e1}) is the rate of growth,
$\mu,$ of the mean population size, i.e.,
\begin{equation}
\log\mu:=\underset{t\rightarrow\infty}{\lim}\log\mathbb{E}\left\Vert
\boldsymbol{z}(t)\right\Vert /t. \label{logmu}%
\end{equation}

It is well known that, as an immediate consequence of Jensen's inequality
\citep{billingsley2008probability}, one has
\begin{equation}
\lambda_{\mathrm{S}}\leq\mu, \label{desiglogmu}%
\end{equation}
and so $\mu$ is an upper bound for the SGR. Note that one can have apparently
paradoxical situations in which $\log\mu>0$ and so the mean population size
grows to infinity, whereas $\log\lambda_{\mathrm{S}}<0,$ and so the population
goes extinct with probability one. This illustrates the sufficient yet
unnecessary essence of condition $\log\mu<0$ for population extinction.

The literature shows \citep{cohen1979long} that $\mu$ coincides with the
dominant eigenvalue of a certain nonnegative matrix $\boldsymbol{D}$ that, in
the general Markovian case, depends on matrices $\boldsymbol{A}_{\eta}$ and on
the one-step transition probabilities of the chain $\tau_{t}$. If, for the
sake of simplicity, we restrict our attention to the case in which $\tau_{t}$
is an independent and identically distributed (IID) process, $\boldsymbol{D}$
coincides with the mean environmental matrix, i.e.,
\begin{equation}
\boldsymbol{\bar{A}}:=\sum_{\eta\in\mathcal{I}}\pi_{\eta}\boldsymbol{A}_{\eta
}. \label{meanmatrix}%
\end{equation}

It is easy to check that when there is no environmental variability, matrix
$\boldsymbol{D}$ coincides with the deterministic matrix $\boldsymbol{A}$ and
so $\mu$ verifies property P.

In the case of system (\ref{Leslie}) and IID temporal variation, one has%
\[
\log\mu=\log\rho(\mathbf{\bar{A}})=\log\frac{\bar{f}+\left(  \bar{f}^{2}%
+4\bar{F}\bar{s}\right)  ^{1/2}}{2}.
\]
When imposing the sufficient condition for extinction $\mu<1,$ the resulting
analytical expression can be simplified by using the concept of \textit{net
reproductive value} in matrix models \citep{cushing1994net}. Indeed if
$R_{0}=\bar{f}+\bar{s}\bar{F}$ denotes the net reproductive value
corresponding to matrix $\boldsymbol{\bar{A}}$, the condition $\mu<1$ is
equivalent to $R_{0}<1$, i.e., $\bar{f}+\bar{s}\bar{F}<1$ which is, therefore,
a sufficient condition for the extinction of the population.

\noindent\textit{\textbf{Cohen's bounds}}. For each $\eta\in\mathcal{I}$,
define the column sums of $\boldsymbol{A}_{\eta}$
\begin{equation}
\varphi_{\eta}^{j}:=\overset{n}{\underset{i=1}{\sum}}A_{\eta}^{ij}%
,\ j=1,...,n\text{, }\eta\in\mathcal{I}, \label{sumas}%
\end{equation}
and let $m_{\eta}$ and $M_{\eta}$ be their respective minimum and maximum,
i.e.,
\[
m_{\eta}=\underset{j=1,...,n}{\min}\varphi_{\eta}^{j},\ \ \ M_{\eta}%
=\underset{j=1,...,n}{\max}\varphi_{\eta}^{j},\ \eta\in\mathcal{I}.
\]
Then \citep{cohen1979long} one has $c_{I}\leq\log\lambda_{\mathrm{S}}\leq
C_{I}$ where
\begin{equation}
c_{I}:=\underset{\eta\in\mathcal{I}}{\sum}\pi_{\eta}\log m_{\eta}%
,\ \ C_{I}:=\underset{\eta\in\mathcal{I}}{\sum}\pi_{\eta}\log M_{\eta}.
\label{CotasCohen}%
\end{equation}

In \citep{tuljapurkar1980population} certain numerical simulations are carried
out that show that, at least in the cases covered in the study, bounds
(\ref{CotasCohen}) do not work well. Bounds (\ref{CotasCohen}) are tight when
for each $\eta\in\mathcal{I}$ all the columns of $A_{\eta}$ sum approximately
to the same amount, and the contrary happens when the maximum and the minimum
of the column sums are very different. Indeed, it is easy to check that they
are exact if and only if $\boldsymbol{A}_{\eta}=\alpha_{\eta}\boldsymbol{P}$
for each $\eta\in\mathcal{I}$, where $\alpha_{\eta}\in\mathbb{R}$ and
$\boldsymbol{P}$ is a column-stochastic matrix.

Another problem of bounds (\ref{CotasCohen}) is that they do not verify
property P, i.e., they do not coincide with the SGR when there is no
environmental variation.

In the case of system (\ref{Leslie}), one has%
\begin{equation}
c_{I}:=\underset{\eta\in\mathcal{I}}{\sum}\pi_{\eta}\log\min\left\{  f_{\eta
}+s_{\eta},F_{\eta}\right\}  ,\ C_{I}:=\underset{\eta\in\mathcal{I}}{\sum}%
\pi_{\eta}\log\max\left\{  f_{\eta}+s_{\eta},F_{\eta}\right\}  .
\label{cotasI}%
\end{equation}

\noindent\textit{\textbf{Bounds through a max-min approach}}. It is immediate
to show that if $\boldsymbol{\tilde{A}},\boldsymbol{\hat{A}}\in\mathbb{R}%
^{n\times n}$ are non-negative matrices such that
\[
0\leq\boldsymbol{\tilde{A}}\leq\boldsymbol{A}_{\eta}\leq\boldsymbol{\hat{A}%
},\ \eta\in\mathcal{I},
\]
then $c_{\min}\leq\log\lambda_{\mathrm{S}}\leq C_{\max}$ where
\begin{equation}
c_{\min}:=\log\rho(\boldsymbol{\tilde{A}}),\ C_{\max}:=\log\rho
(\boldsymbol{\hat{A}}) \label{CotasMaxMin}%
\end{equation}
Clearly, these bounds are very easy to obtain and verify property P. A serious
drawback is that they are independent of the equilibrium environmental
distribution $\boldsymbol{\pi}$ and so they do not weight the relative
contribution of the different environments. For instance, when there is only
two environments and $\boldsymbol{A}_{1}<\boldsymbol{A}_{2}$, it is easy to
see that $c_{\min}$ is very poor whenever $\pi_{1}$ is close to zero.

In the case of system (\ref{Leslie}), one has%
\begin{equation}
\boldsymbol{\tilde{A}}=\left(
\begin{array}
[c]{cc}%
\underset{^{\eta\in\mathcal{I}}}{\min}f_{\eta} & \underset{^{\eta
\in\mathcal{I}}}{\min}F_{\eta}\\
\underset{^{\eta\in\mathcal{I}}}{\min}s_{\eta} & 0
\end{array}
\right)  ,\ \boldsymbol{\hat{A}}=\left(
\begin{array}
[c]{cc}%
\underset{^{\eta\in\mathcal{I}}}{\max}f_{\eta} & \underset{^{\eta
\in\mathcal{I}}}{\max}F_{\eta}\\
\underset{^{\eta\in\mathcal{I}}}{\max}s_{\eta} & 0
\end{array}
\right)  . \label{cotasII}%
\end{equation}

\section{A perturbation approach to obtaining bounds for the
SGR\label{sec:pert}}

Let us first introduce a new result that relates the SGR of any two stochastic
systems of the kind (\ref{e1}). Although more general in nature, this result
will be particularized to certain situations to obtain bounds for the SGR.

Consider the following stochastic system:
\begin{equation}
\boldsymbol{z}^{\prime}(t+1)=\boldsymbol{B}_{\tau_{t+1}}\boldsymbol{z}%
^{\prime}(t), \label{e1bis}%
\end{equation}
where we assume that $\left\{  \boldsymbol{B}_{\eta},\eta\in\mathcal{I}%
\right\}  $ is an ergodic set of $n\times n$ non-negative matrices and,
further, that their incidence matrices coincide with those of system
(\ref{e1}), that is, for all $\eta\in\mathcal{I}$ one has $B_{\eta}^{ij}=0$ if
and only if $A_{\eta}^{ij}=0$.

Under the above conditions there exists an SGR, that we denote $\lambda
_{\mathrm{S}}^{\prime},$ for (\ref{e1bis}). Now we want to relate $\lambda
_{S}$ and $\lambda_{\mathrm{S}}^{\prime}$. In order to do so, let us first
define matrices $\boldsymbol{W}_{\eta}=\left[  W_{\eta}^{ij}\right]  _{1\leq
i,j\leq n}$, $\eta\in\mathcal{I}$ through
\begin{equation}
W_{\eta}^{ij}:=\left\{
\begin{tabular}
[c]{c}%
$\frac{A_{\eta}^{ij}-B_{\eta}^{ij}}{B_{\eta}^{ij}}\text{ if }B_{\eta}^{ij}%
\neq0$\\
$0\text{ if }B_{\eta}^{ij}=0$%
\end{tabular}
\ \ .\right.  \label{pert02}%
\end{equation}
Note that $W_{\eta}^{ij}$ can be thought of as the relative amount of
perturbation of $A_{\eta}^{ij}$ with respect to $B_{\eta}^{ij}$. Let us now
set
\begin{equation}
W_{\eta}^{m}:=\min_{i,j=1,...,n}W_{\eta}^{ij},\ \ W_{\eta}^{M}:=\max
_{i,j=1,...,n}W_{\eta}^{ij}. \label{pert03}%
\end{equation}
Using the nonnegativity of $\boldsymbol{A}_{\eta}$ and $\boldsymbol{B}_{\eta}$
and the fact that the incidence matrix of $\boldsymbol{B}_{\eta}$ coincides
with that of $\mathbf{A}_{\eta}$, we check immediately that
\begin{equation}
W_{\eta}^{M}\geq W_{\eta}^{m}>-1. \label{pert04}%
\end{equation}

Then (see \ref{sec:proofs}) one has%
\begin{equation}
\log\lambda_{\mathrm{S}}^{\prime}+\underset{\eta\in\mathcal{I}}{\sum}\pi
_{\eta}\log\left(  1+W_{\eta}^{m}\right)  \leq\log\lambda_{S}\leq\log
\lambda_{\mathrm{S}}^{\prime}+\underset{\eta\in\mathcal{I}}{\sum}\pi_{\eta
}\log\left(  1+W_{\eta}^{M}\right)  . \label{CotasPertBis}%
\end{equation}

The previous result is more general, but in this work it is used in the case
in which system (\ref{e1bis}) is deterministic, i.e., $\boldsymbol{B}_{\eta
}=\boldsymbol{B}$ for all $\eta\in\mathcal{I}$ where $\boldsymbol{B}$ is a
certain matrix, so that $\lambda_{\mathrm{S}}^{\prime}=\rho(\boldsymbol{B})$.
Therefore in that case one can write $c_{\text{pert}}\leq\log\lambda_{S}\leq
C_{\text{pert}}$ where%
\begin{equation}
c_{\text{pert}}:=\log\rho(\boldsymbol{B})+\underset{\eta\in\mathcal{I}}{\sum
}\pi_{\eta}\log\left(  1+W_{\eta}^{m}\right)  ,\ C_{\text{pert}}:=\log
\rho(\boldsymbol{B})+\underset{\eta\in\mathcal{I}}{\sum}\pi_{\eta}\log\left(
1+W_{\eta}^{M}\right)  . \label{CotasPert}%
\end{equation}

Bounds $c_{\text{pert}}$ and $C_{\text{pert}}$ have been deduced through a
perturbation approach and certainly work better when, for each $\eta
\in\mathcal{I},$ matrix $\boldsymbol{A}_{\eta}$ is a perturbation of
$\boldsymbol{B}$. However, the bounds are valid even when that is not the
case. We stress that these bounds are valid in systems of general
dimensionality, i.e., $n$ is arbitrary.

Note that when there is no stochasticity $W_{\eta}^{m}=W_{\eta}^{M}=0$ for all
$\eta\in\mathcal{I},$ and so the bounds in (\ref{CotasPert}) verify property P
as desired.

Obviously different choices of matrix $\boldsymbol{B}$ will lead to different
bounds. We propose three possibilities:

\noindent1. $\boldsymbol{B}$ \textit{\textbf{equals the average environmental
matrix.}} We take
\[
\boldsymbol{B}=\boldsymbol{\bar{A}}:=\mathbb{E}_{\pi}(\boldsymbol{A}_{\tau
_{t}})=\underset{\eta\in\mathcal{I}}{\sum}\pi_{\eta}\boldsymbol{A}_{\eta},
\]
i.e., the average environmental matrix with the average taken with respect to
the stationary distribution of $\tau_{t}$. We denote $c_{II}$ and $C_{II}$ the
resulting bounds in (\ref{CotasPert}) for this choice of $\boldsymbol{B}$.

In the case of system (\ref{Leslie}) we have
\begin{align*}
c_{II}  &  =\log\frac{\bar{f}+\left(  \bar{f}^{2}+4\bar{F}\bar{s}\right)
^{1/2}}{2}+\underset{\eta\in\mathcal{I}}{\sum}\pi_{\eta}\log\min\left\{
1,\frac{f_{\eta}}{\bar{f}},\frac{F_{\eta}}{\bar{F}},\frac{s_{\eta}}{\bar{f}%
}\right\}  ,\\
C_{II}  &  =\log\frac{\bar{f}+\left(  \bar{f}^{2}+4\bar{F}\bar{s}\right)
^{1/2}}{2}+\underset{\eta\in\mathcal{I}}{\sum}\pi_{\eta}\log\max\left\{
1,\frac{f_{\eta}}{\bar{f}},\frac{F_{\eta}}{\bar{F}},\frac{s_{\eta}}{\bar{f}%
}\right\}  .
\end{align*}

Note that for system (\ref{Leslie}) and the IID setting, one has $C_{II}%
\geq\log\rho(\boldsymbol{\bar{A}})=\log\mu$ and therefore in this case
$\log\mu$ is always better than $C_{II}$ as an upper bound for $\log
\lambda_{\mathrm{S}}$.

\noindent2. $\boldsymbol{B}$\textit{\textbf{ equals the \textquotedblleft
max-matrix\textquotedblright.}} We define $\boldsymbol{B}=\left[
B^{ij}\right]  _{1\leq i,j\leq n}$ through%
\[
B^{ij}=\underset{\eta\in\mathcal{I}}{\max}A_{\eta}^{ij},
\]
and we denote $c_{III}$ and $C_{III}$ the resulting bounds in (\ref{CotasPert}).

In this case $\boldsymbol{B}=\boldsymbol{\hat{A}}$ in (\ref{CotasMaxMin}) and
$W_{\eta}^{M}\leq0$ for all $\eta\in\mathcal{I}$. Therefore $C_{III}\leq
\log\rho(\boldsymbol{\hat{A}})=C_{\max}$ and so $C_{III}$ is always a better
upper bound than $C_{\max}$.

For system (\ref{Leslie}) one has $W_{\eta}^{M}=0$ for all $\eta\in
\mathcal{I}$ and so%
\begin{gather}
c_{III}=\log\frac{\max\nolimits_{\eta}f_{\eta}+\left(  \left(  \max
\nolimits_{\eta}f_{\eta}\right)  ^{2}+4\left(  \max\nolimits_{\eta}F_{\eta
}\right)  \left(  \max\nolimits_{\eta}s_{\eta}\right)  \right)  ^{1/2}}%
{2}+\label{cotin}\\
+\underset{\eta\in\mathcal{I}}{\sum}\pi_{\eta}\log\left(  1+\min\left\{
\frac{f_{\eta}}{\max\nolimits_{\eta}f_{\eta}},\frac{F_{\eta}}{\max
\nolimits_{\eta}F_{\eta}},\frac{s_{\eta}}{\max\nolimits_{\eta}s_{\eta}%
}\right\}  \right)  ,\nonumber\\
C_{III}=C_{\max}=\log\frac{\max\nolimits_{\eta}f_{\eta}+\left(  \left(
\max\nolimits_{\eta}f_{\eta}\right)  ^{2}+4\left(  \max\nolimits_{\eta}%
F_{\eta}\right)  \left(  \max\nolimits_{\eta}s_{\eta}\right)  \right)  ^{1/2}%
}{2}.\nonumber
\end{gather}

\noindent3. $\boldsymbol{B}$\textit{\textbf{ equals the \textquotedblleft
min-matrix\textquotedblright.}} We put%
\[
B^{ij}=\underset{\eta\in\mathcal{I}}{\min}A_{\eta}^{ij},
\]
and the resulting bounds in (\ref{CotasPert}) are denoted $c_{IV}$ and
$C_{IV}$. Reasoning like in the previous case, one has $c_{IV}\geq\log
\rho(\boldsymbol{\tilde{A}})=c_{\min}$ and therefore $c_{IV}$ is always a
better lower bound than $c_{\min}$. In the case of system (\ref{Leslie}), one
has analogous expressions for $c_{IV}$ and $C_{IV}$ to those of (\ref{cotin}).

We emphasize that the bounds $C_{III}$ and $c_{IV}$ proposed in this work
improve the existing bounds $C_{\max}$ and $c_{\min}$ respectively.

There is no rule of thumb for the best choice of matrix $B,$ even in very
simple situations as when there are only two stages, only one vital rate
varies and there are two environments which are equiprobable. Simulations show
that, in general, the optimal choice of $B$ is somewhere inbetween the
\textquotedblleft min\textquotedblright\ and the \textquotedblleft
max\textquotedblright\ matrices, sometimes close to the former, some to the
latter and some others to the average matrix.

\section{Bounds in two-class Leslie models \label{sec:boundsagestructure}}

The idea is to find bounds for the vector of age structure, i.e., the vector
with the proportion of young and adults, that are valid for large times and
that can be used later on to obtain bounds for the SGR. As a first step, we
generalize some results from \citep{Tulja90} to obtain bounds for the vector
of age structure when all three vital rates (fertility of the young, of the
adults and survival) are stochastic. These bounds on age structure are then
used to obtain a new bound for the SGR.

\subsection{Bounds on age structure for two-class Leslie
models\label{sec:BoundsStr}}

We are interested in finding bounds for age structure valid for large values
of $t$, i.e., we want to find (deterministic) vectors $\boldsymbol{l}=\left(
l^{1},...,l^{n}\right)  ,\boldsymbol{u}=\left(  u^{1},...,u^{n}\right)
\in\mathbb{R}_{+}^{n},$ such that for large enough $t$ one can write
\begin{equation}
\boldsymbol{l}\leq\frac{\boldsymbol{z}(t)}{\left\Vert \boldsymbol{z}%
(t)\right\Vert }\leq\boldsymbol{u} \label{wer}%
\end{equation}
These bounds on age structure are relevant in this work since, as it is shown
later on, they can be used to obtain new bounds for the SGR.

It is not possible to obtain bounds of the kind (\ref{wer}) for system
(\ref{e1}) in the general case. In \citep{Tulja90} certain bounds are obtained
in the particular case of Leslie models, although when there is an arbitrary
number of age classes the bounds can not be obtained through closed
expressions but must be found through an iterative process. Here we
concentrate in the case of system (\ref{Leslie}) with $n=2$ age classes, for
which closed form expressions can be found.

First we obtain bounds for the quotient $z^{2}(t)/z^{1}(t)$ in the form
\begin{equation}
\delta\leq\frac{z^{2}(t)}{z^{1}(t)}\leq\kappa. \label{CotasEst}%
\end{equation}
Once this is accomplished, it is immediate to conclude that the vector of
population structure can be bounded as follows
\begin{equation}
\left(  \frac{\delta}{1+\kappa},\frac{1}{1+\kappa}\right)  =:\boldsymbol{l}%
\leq\frac{\boldsymbol{z}(t)}{\left\Vert \boldsymbol{z}(t)\right\Vert }%
\leq\boldsymbol{u}:=\left(  \frac{\kappa}{1+\delta},\frac{1}{1+\delta}\right)
. \label{para}%
\end{equation}

In \citep{Tulja90} bounds of the kind (\ref{CotasEst}) are obtained in the
case in which only the survival of the young is stochastic. Here we show that,
following a similar reasoning to that of \citep{Tulja90}, the results can be
extended to the case in which the three vital rates are stochastic.
Specifically, let us define
\[
\varepsilon_{M}:=\underset{\eta\in\mathcal{I}}{\max}\frac{F_{\eta}}{f_{\eta}%
},\ \varepsilon_{m}:=\underset{\eta\in\mathcal{I}}{\min}\frac{F_{\eta}%
}{f_{\eta}},\ \gamma_{M}:=\underset{\eta\in\mathcal{I}}{\max}\frac{s_{\eta}%
}{f_{\eta}},\ \gamma_{m}:=\underset{\eta\in\mathcal{I}}{\min}\frac{s_{\eta}%
}{f_{\eta}}.
\]
Then it can be shown (see \ref{sec:proofs} for the proof) that equation
(\ref{CotasEst}) holds with
\begin{align}
\delta &  =\frac{2\gamma_{m}}{1-\varepsilon_{m}\gamma_{m}+\varepsilon
_{M}\gamma_{M}+\sqrt{1+(\varepsilon_{M}\gamma_{M}-\varepsilon_{m}\gamma
_{m})^{2}+2(\varepsilon_{M}\gamma_{M}+\varepsilon_{m}\gamma_{m})}%
},\label{cotitas}\\
\kappa &  =\frac{2\gamma_{M}}{1-\varepsilon_{M}\gamma_{M}+\varepsilon
_{m}\gamma_{m}+\sqrt{1+(\varepsilon_{m}\gamma_{m}-\varepsilon_{M}\gamma
_{M})^{2}+2(\varepsilon_{m}\gamma_{m}-\varepsilon_{M}\gamma_{M})}}.\nonumber
\end{align}

Note that, in the deterministic setting, i.e., $\boldsymbol{A}_{\eta
}=\boldsymbol{A}$ for all $\eta\in\mathcal{I}$, bounds $\delta$ and $\kappa$
coincide with the quotient $v^{2}/v^{1}$ where $\boldsymbol{v}=(v^{1},v^{2})$
is a right eigenvector of matrix $\boldsymbol{A}$ associated to its dominant
eigenvalue, and so bounds (\ref{para}) are exact.

\subsection{Bounds for the SGR using info on age
structure\label{sec:CotasCohenEstr}}

Consider system (\ref{e1}) and assume that (\ref{wer}) holds for certain known
vectors $\boldsymbol{l}$ and $\boldsymbol{u}$. Then it follows (see
\ref{sec:proofs} for the proof) that $c_{V}\leq\log\lambda_{\mathrm{S}}\leq
C_{V}$ where
\begin{equation}
c_{V}:=\underset{\eta\in\mathcal{I}}{\sum}\pi_{\eta}\log\left(  \overset
{n}{\underset{j=1}{\sum}}l_{\eta}^{j}\varphi_{\eta}^{j}\right)  ,\ C_{V}%
:=\underset{\eta\in\mathcal{I}}{\sum}\pi_{\eta}\log\left(  \overset
{n}{\underset{j=1}{\sum}}u^{j}\varphi_{\eta}^{j}\right)  ,
\label{CotasCohenEstr}%
\end{equation}
and the $\varphi_{\eta}^{j}$ are given by (\ref{sumas}).

Note that in the deterministic case one has $\boldsymbol{l}=\boldsymbol{u}$
and so bounds (\ref{CotasCohenEstr}) are exact.

In principle these bounds can be used for any system of the kind (\ref{e1}),
but as mentioned before, in practice bounds of the kind (\ref{wer}) are not at
our disposal. In the particular case of model (\ref{Leslie}), by making use of
(\ref{para}) one obtains%
\[
c_{V}:=\underset{\eta\in\mathcal{I}}{\sum}\pi_{\eta}\log\left(  \frac{f_{\eta
}+s_{\eta}+F_{\eta}\delta}{1+\kappa}\right)  ,\ C_{V}:=\underset{\eta
\in\mathcal{I}}{\sum}\pi_{\eta}\log\left(  \frac{f_{\eta}+s_{\eta}+F_{\eta
}\kappa}{1+\delta}\right)  ,
\]
where $\delta$ and $\kappa$ are given by (\ref{cotitas}).

\section{Analysis and comparison of the different bounds\label{sec:comparison}%
}

In this section we carry out a comparison of the different bounds which have
been previously introduced. Before we proceed, it should be noted that,
excluding $\mu$, the rest of the bounds considered in this work depend on the
environmental chain $\tau_{t}$ only through its equilibrium distribution
$\boldsymbol{\pi}$ and are therefore independent of the serial correlation
between different environments. On the other hand, we stress that previous
studies \citep{tuljapurkar1982population} show that the SGR depends strongly
on the stationary distribution $\boldsymbol{\pi}$ but rather weakly on the
serial correlation between environments. Taking these two facts into
consideration, in order to carry out the comparison between the different
bounds, we restrict our attention to the IID setting, i.e., there is no serial
correlation between environments. More specifically, we consider the two class
Leslie model (\ref{Leslie}) with two different environments. Note that the
environmental process $\tau_{t}$ is completely defined through the value of
$\pi_{1}\in(0,1)$.

In particular, we want to find out which upper and which lower bounds wins for
each value of the parameter values, i.e., which lower bound is the largest and
which upper bound is the lowest. Even in this simple setting, the large number
of parameters involved makes an analytical study of the comparison unfeasible,
so resort to numerical simulations. Around $2\ast10^{6}$ numerical simulations
have been performed in which the parameters $\pi_{1}$, $f_{1}$, $f_{2}$,
$F_{1}$, $F_{2}$, $s_{1}$ and $s_{2}$ take different values. Specifically,
$\pi_{1}\in\{0.1,\ 0.3,\ 0.5,\ 0.7,\ 0.9\}$, and the vital rates are taken
from the following intervals: $f_{1},f_{2}\in\lbrack0.1,3]$, $F_{1},F_{2}%
\in\lbrack0.1,5]$, $s_{1},s_{2}\in\lbrack0.1,0.9]$. The distance between
consecutive values of each parameter is 0.2.

\begin{table}[pt]
\centering
\begin{tabular}
[c]{|c|c|c|c|c|c|}\hline
$C_{I}$ & $C_{II}$ & $C_{III}$ & $C_{IV}$ & $C_{V}$ & $\log\mu$\\\hline
2.57\% & 0\% & 2.78\% & 0\% & 2.44\% & 92.20\%\\\hline
\end{tabular}
\caption{Percentage of cases in which each upper bound wins.}%
\label{tabla2}%
\end{table}

\begin{table}[pt]
\centering
\begin{tabular}
[c]{|c|c|c|c|c|}\hline
$c_{I}$ & $c_{II}$ & $c_{III}$ & $c_{IV}$ & $c_{V}$\\\hline
36.44\% & 8.62\% & 7.90\% & 18.75\% & 28.30\%\\\hline
\end{tabular}
\caption{Percentage of cases in which each lower bound wins.}%
\label{tabla1}%
\end{table}

Tables \ref{tabla2} and \ref{tabla1} show the percentage of cases in which
each upper bound and each lower wins. Regarding the upper bounds we can draw
the following conclusions:

\begin{itemize}
\item $\mu$ is the best upper bound in the vast majority of cases. This,
combined with the fact that its analytical expression is relatively simple,
makes it the best choice of upper bound when one does not want to check which
upper bound is the best.

\item Bounds $C_{II}$ and $C_{IV}$ never win.

\item Bounds $C_{I}$, $C_{III}$ and $C_{V}$ win in approximately the same
percentage of cases.
\end{itemize}

Regarding the lower bounds, one can see that $c_{I}$, $c_{IV}$ and $c_{V}$ win
in a large percentage of cases, but so are $c_{II}$ and $c_{III}$ for certain
instances. In Section \ref{sec:delta} we illustrate this fact with specific examples.

Simulations have also shown that when the two environments are more or less
equiprobable, i.e., $\pi_{1}\approx0.5$, then $c_{I}$ and $c_{V}$ win in most
cases (45.16\% and 35.75\% respectively for $\pi_{1}=0.5$) and $c_{II}$ never
win, but when one environment is much more frequent, for example for $\pi
_{1}=0.1$, $c_{III}$ wins in 10\% of cases and, in the remaining 90\%, the
rest of the bounds win in approximately in the same percentage of cases.

The second question we want to analyze is how large is the relative error
incurred when using the best upper or best lower bound for $\lambda
_{\mathrm{S}}$. Note that we work with $\lambda_{\mathrm{S}},$ and not with
$\log\lambda_{\mathrm{S}},$ to avoid having to divide by numbers close to zero
when computing the relative error in the cases for which $\lambda_{\mathrm{S}%
}$ is close to $1$.

It is reasonable to expect that this error depends strongly on the amount of
environmental variation. Therefore, a set of simulations similar to the ones
described above has been carried out. We have computed the relative error
$100\ast\left\vert \lambda_{\mathrm{S}}-C\right\vert /\lambda_{\mathrm{S}},$
where $C$ is the best (upper or lower) bound for $\lambda_{\mathrm{S}}$, and
have compared it to the relative amount of environmental variation, measured
as
\[
100\ast\frac{\left\vert f_{2}-f_{1}\right\vert +\left\vert s_{2}%
-s_{1}\right\vert +\left\vert F_{2}-F_{1}\right\vert }{\max\left\{
f_{1},f_{2}\right\}  +\max\left\{  s_{1},s_{2}\right\}  +\max\left\{
F_{1},F_{2}\right\}  }.
\]

The SGR can not be obtained analytically and, therefore, has been approximated
by means of numerical simulations. Specifically, we have used the approach
described in \cite{tuljapurkar1997stochastic} and, in order to keep the
computational cost manageable and still obtain a good approximation, be have
used $N=$500 samples each of them with a length of $T=$600 time steps. In
order to estimate the sampling error for these values of $N$ and $T$ we have
compared the results for these values with the ones corresponding to very
large values of these parameters, specifically $N=2000$ and $T=8000,$ for
which we can expect a very accurate estimate of $\lambda_{\mathrm{S}}$ to be
obtained. Thus, we have computed $\lambda_{\mathrm{S}}$ in 100 different cases
using both procedures, and have obtained that the maximum discrepancy between
them is around 0.3\%, whereas the mean discrepancy is around 0.1\%. Therefore,
we can consider that the sampling error has almost no effect in the estimation
of the relative error for the best upper and lower bounds.

\begin{figure}[h]
\centering
\begin{subfigure}[t]{0.45\textwidth}
\centering
\includegraphics[width=\linewidth]{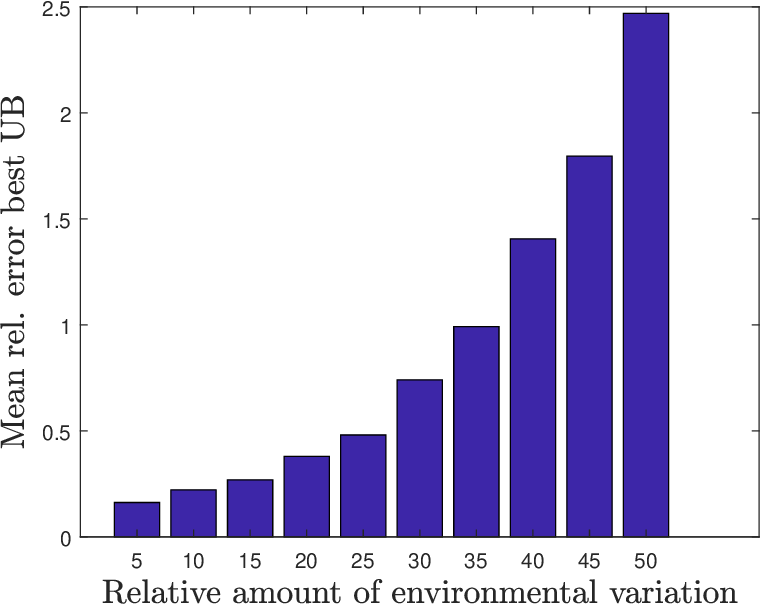}
\end{subfigure}
\qquad\begin{subfigure}[t]{0.45\textwidth}
\centering
\includegraphics[width=\linewidth]{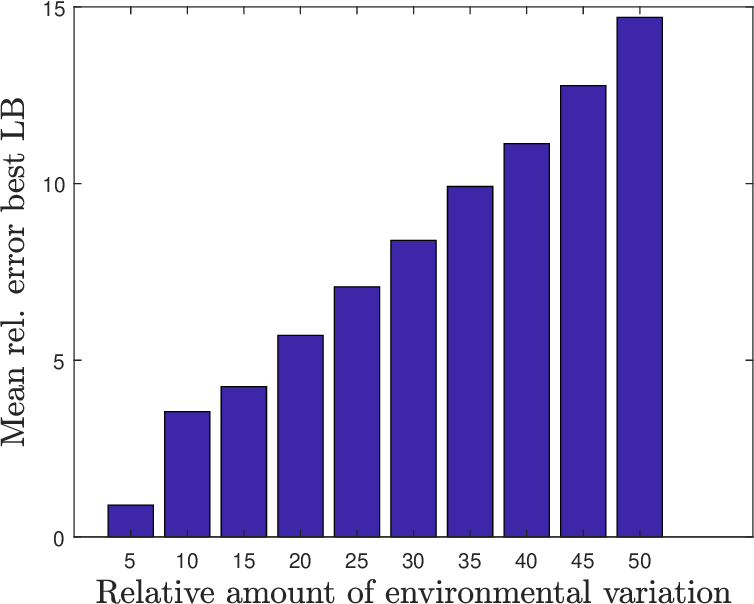}
\end{subfigure}
\par
\rule{0ex}{2ex}
\par
\centering
\begin{subfigure}[t]{0.45\textwidth}
\centering
\includegraphics[width=\linewidth]{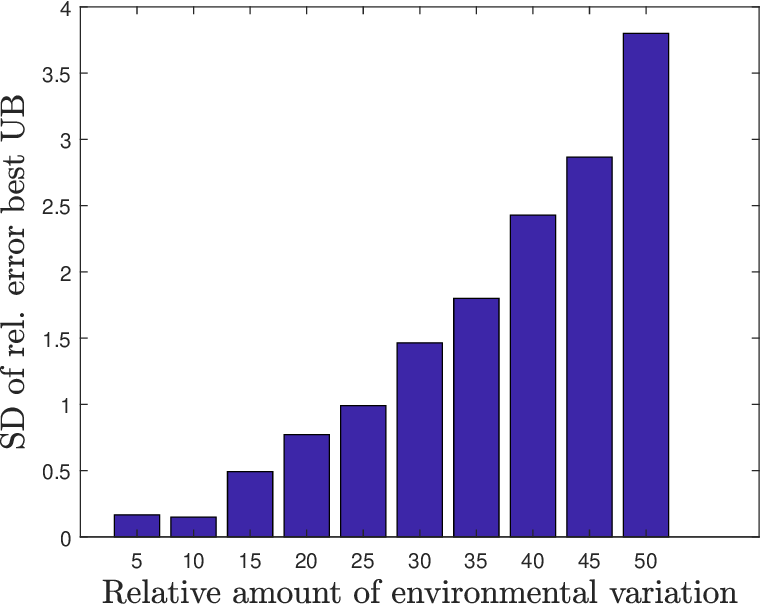}
\end{subfigure}
\qquad\begin{subfigure}[t]{0.45\textwidth}
\centering
\includegraphics[width=\linewidth]{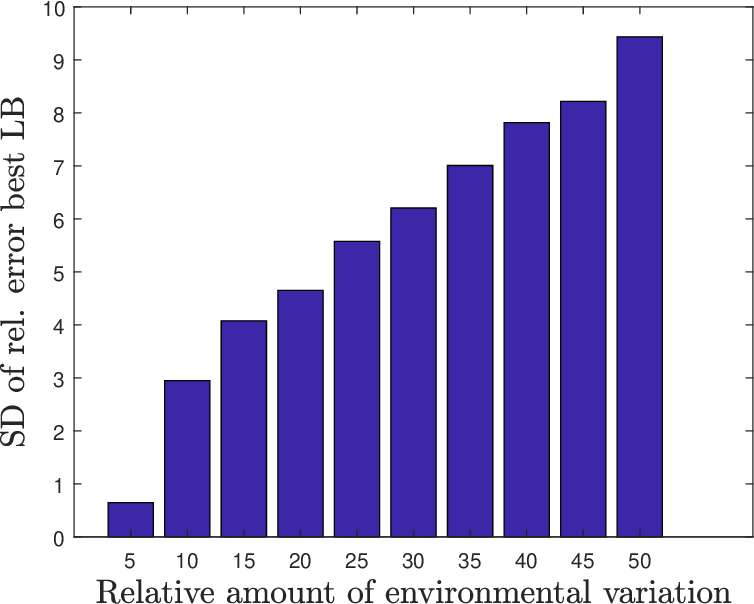}
\end{subfigure}
\caption{Mean and standard deviation (SD) of relative error (\%) for the best
upper bound (UB) and the best lower bound (LB) of $\lambda_{S}$ as functions
of relative environmental variation (\%)}%
\label{Fig2}%
\end{figure}

The results are shown in Figure \ref{Fig2}. They suggest that both the best
upper bound and the best lower bound are quite good approximations to the SGR
for a wide range of values of environmental variation. They also show that the
mean error for the best upper bound is roughly half that corresponding to the
best lower bound. Note the comparatively large standard deviation in the
distribution of the error corresponding to both the best upper and best lower
bound, specially in the latter case.

We conclude this section with some comments about Tuljapurkar's second order
approximation to the SGR $\lambda_{T}$ \citep{tuljapurkar1982population}
which, as pointed out before, is a popular tool for ecologists. In the
previous simulations we have analyzed some properties of $\lambda_{T}$ and
have reached interesting conclusions. In the first place, in around 95\% of
cases $\lambda_{T}$ underestimates the SGR. Secondly, in around 99\% of cases
$\lambda_{T}$ is larger than the best lower bound. These two points suggest
that, as a rule of thumb that works in most situations, we could consider
$\lambda_{T}$ as a \textquotedblleft practical\ lower bound\textquotedblright%
\ for the SGR (although from a rigorous point of view it is not) that, in many
cases, performs better that the lower bounds that have been studied in this work.

\section{A study of growth/extinction in a stochastic Leslie
model\label{sec:delta}}

As an illustration of a practical use of the bounds presented in this work,
let us consider the Leslie model with two age classes (\ref{Leslie}) in the
following setting: there are two possible environments, i.e., $\mathcal{I}%
=\left\{  1,2\right\}  $. Environment $1$ is rich in the sense that
$\rho(\boldsymbol{A}_{1})>1$ whereas environment $2$ is poorer in the sense
that the adults fertility is reduced by an amount $\Delta>0$, i.e.,%
\[
\boldsymbol{A}_{1}:=\left(
\begin{array}
[c]{cc}%
f & F\\
s & 0
\end{array}
\right)  ,\ \boldsymbol{A}_{2}:=\left(
\begin{array}
[c]{cc}%
f & F-\Delta\\
s & 0
\end{array}
\right)  .
\]

We consider fixed values of $f>0$, $F>0$ and $s\in(0,1],$ with the restriction
that $\rho(A_{1})>1$ or alternatively, using the \textit{net reproductive
value}, that
\begin{equation}
f+sF>1. \label{crec}%
\end{equation}
Therefore, if subjected to environment 1 alone, the population grows to
infinity. Let us consider the IID setting, let $\pi_{1}\in(0,1)$ be fixed and
let $\boldsymbol{\pi}=\left(  \pi_{1},1-\pi_{1}\right)  $ be the probability
distribution of environmental states.

Clearly, when $\Delta=0,$ we only have environment 1 and the population grows
to infinity. Also, as $\Delta$ increases the population fitness decreases.
Therefore, we want to answer the following question: what is the maximum value
of $\Delta$, that we denote $\Delta_{\text{max}}$, such that $\lambda
_{\mathrm{S}}>1$ and so the population does not go extinct? Since the SGR can
not be derived analytically, the problem can not be solved directly. Our
strategy will be to consider one lower bound for the SGR, say $c<\log
\lambda_{\mathrm{S}}$, and to obtain the maximum value $\Delta_{c}$ such that
$c(\Delta)>1$ whenever $0<\Delta<\Delta_{c}$. Consequently, if the last
inequality holds we know that the population grows to infinity. Note that
$\Delta_{c}\leq\Delta_{\text{max}}$, i.e., with this approach we obtain
sufficient (but not necessary) conditions on $\Delta$ for population growth.
Whenever possible, we want to find closed form expressions for $\Delta_{c}$.
In what follows we will\ complete this process in turn for the different lower
bounds for the SGR.

\noindent\textit{\textbf{Use of Cohen's bound }}$c_{I}$\textit{ }%
(\ref{cotasI}). We distinguish two situations:

\begin{itemize}
\item Case a. $F>1$. In that case it is easy to show, using (\ref{crec}), that
$f+s<F$. Therefore from (\ref{cotasI}) one has $c_{I}=\log\left(  \left(
f+s\right)  ^{\pi_{1}}\left(  F-\Delta\right)  ^{1-\pi_{1}}\right)  $ and
therefore $c_{I}>0$ if and only if
\[
\Delta<\Delta_{IA}:=F-\left(  f+s\right)  ^{\frac{\pi_{1}}{\pi_{1}-1}}.
\]

\item Case b. $F<1$. Then $f+s>F$ and $c_{I}=\log\left(  F^{\pi_{1}}\left(
F-\Delta\right)  ^{1-\pi_{1}}\right)  ,$ so $c_{I}>0$ if and only if
\[
\Delta<\Delta_{IB}:=F\left(  1-F^{\frac{1}{\pi_{1}-1}}\right)  .
\]

\end{itemize}

\noindent\textit{\textbf{Use of bound} }$c_{III}$. From (\ref{cotin}) it
follows that%
\[
c_{III}=\log\rho(\boldsymbol{A}_{1})+\left(  1-\pi_{1}\right)  \log\left(
1-\frac{\Delta}{F}\right)  ,
\]
and so $c_{III}>0$ if and only if
\[
\Delta<\Delta_{III}:=F\left(  1-\rho(\boldsymbol{A}_{1})^{\frac{1}{\pi_{1}-1}%
}\right)  .
\]
\smallskip

\noindent\textit{\textbf{Use of bound} }$c_{IV}$. In this case $c_{IV}%
=c_{\min}=\rho(\boldsymbol{A}_{2}),$ and so $c_{IV}>0$ if and only if
$\rho(\boldsymbol{A}_{2})>1$, i.e., if and only if $f+s(F-\Delta)>1$ or,
equivalently,%
\[
\Delta<\Delta_{IV}:=F-\frac{1-f}{s}.
\]

One can proceed similarly by working with bound $c_{II}$. A closed form
expression of the maximum value of $\Delta$ is not possible using $c_{V}$ due
to the complexity of the equation involved.

\begin{figure}[h]
\centering
\begin{subfigure}[t]{0.70\textwidth}
\centering
\includegraphics[width=\linewidth]{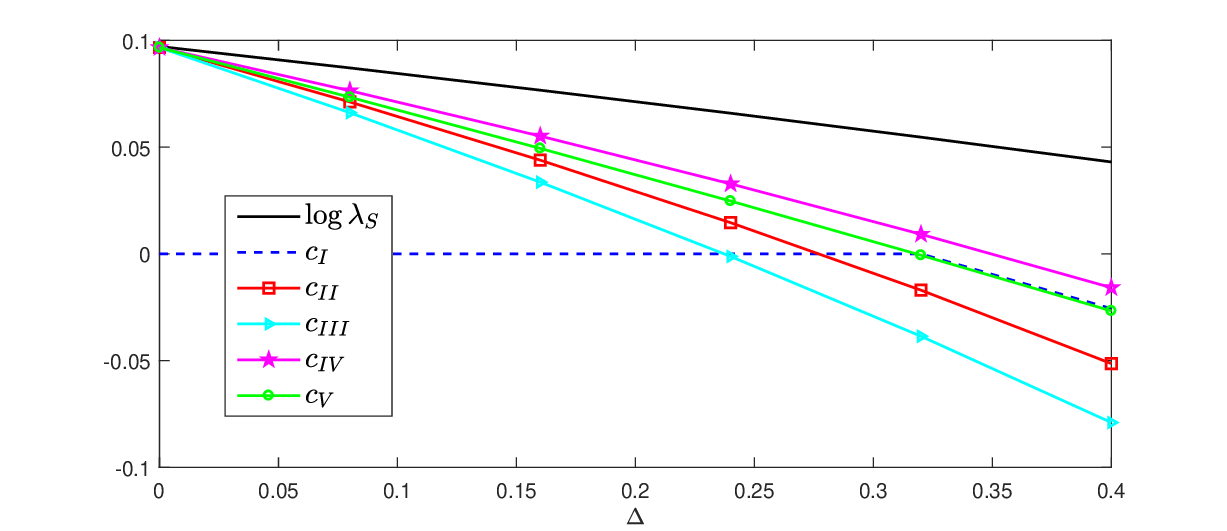}
\caption{\textbf{Case A}. $\pi_1=0.5$, $F=1.35$, $f=0.55$, $s=0.45.$}\label{fig1.1}
\end{subfigure}
\par
\rule{0ex}{2ex}
\par
\begin{subfigure}[t]{0.70\textwidth}
\centering
\includegraphics[width=\linewidth]{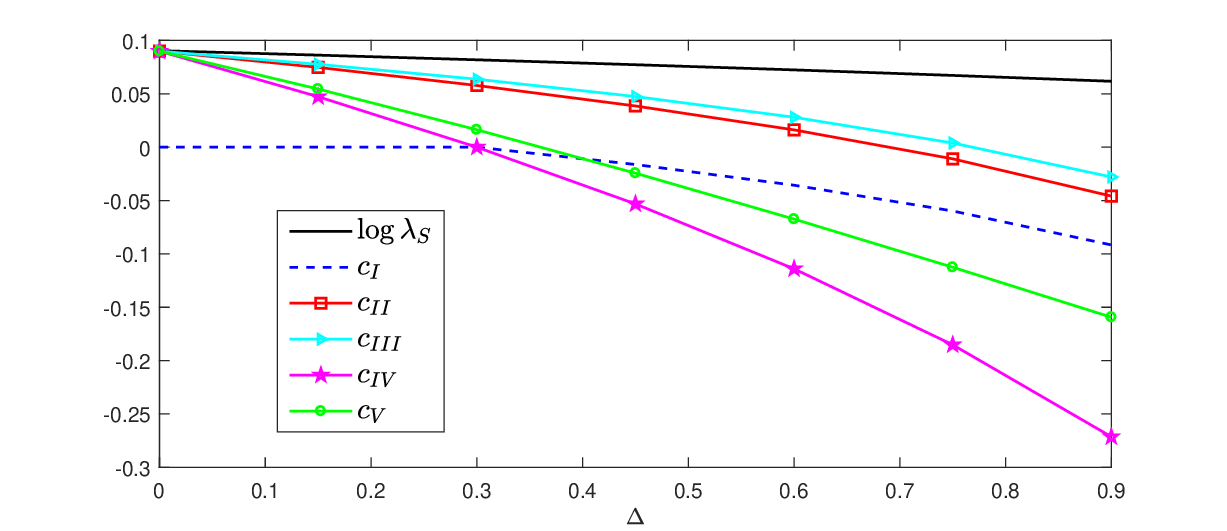}
\caption{\textbf{Case B}. $\pi_1=0.9$, $F=1.3$, $f=0.5$, $s=0.5.$} \label{fig1.2}
\end{subfigure}
\par
\rule{0ex}{2ex}
\par
\centering
\begin{subfigure}[t]{0.70\textwidth}
\centering
\includegraphics[width=\linewidth]{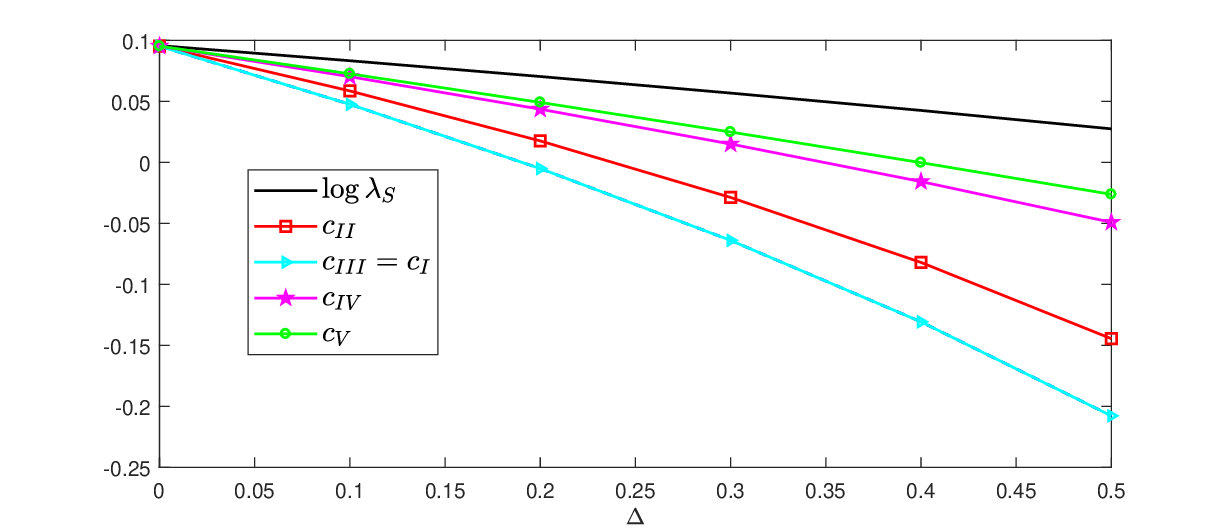}
\caption{\textbf{Case C}. $\pi_1=0.5$, $F=1.1$, $f=0.7$, $s=0.4.$} \label{fig1.3}
\end{subfigure}
\caption{Values of $\log\lambda_{S}$ and of the lower bounds as a function of
$\Delta$ in three different configurations.}%
\label{fig1}%
\end{figure}

The question now is to find out, for each value of the model parameters, which
lower bound wins, i.e., which of the different lower bounds provides the
largest value of $\Delta$ for which one can guarantee non extinction. In order
to address this issue, several numerical simulations have been performed that
show that, even in this particular setting in which only one parameter varies
from environment 1 to environment 2, there is no overall best lower bound, but
different lower bounds can win depending on parameters values. Figure
\ref{fig1} shows the values of $\log\lambda_{\mathrm{S}}$ (obtained through
numerical simulations) and of the lower bounds as functions of $\Delta$ for
three different configurations, that we denote Case A, Case\ B and Case C. In
Case A, $c_{IV}$ wins and the worst bounds are $c_{I}$ for low values of
$\Delta$ and $c_{III}$ for larger values. The situation changes dramatically
in Case B, for $c_{III}$ becomes the winner, whereas the worst bounds are now
$c_{I}$ for low values of $\Delta$ and $c_{IV}$ for larger values. In Case C
bound $c_{V}$ wins whereas $c_{I}$ and $c_{III}$ (in this case these two
bounds coincide) are the worst.

We argued in Section \ref{sec:comparison} that in most practical cases we can
consider $\lambda_{T}$ is a \textquotedblleft practical\textquotedblright%
\ lower bound for $\lambda_{\mathrm{S}}$. Therefore, we could try, in
principle, to find the maximum value of $\Delta$ proceeding as above. However
the resulting expression for $\lambda_{T}$ as a function of $\Delta$ is very
complicated and so it is not possible to find a closed-form expression for the
maximum value of $\Delta$ using $\lambda_{T}$. Numerically, one can observe
that in the three cases A, B and C, $\lambda_{T}$ is lower than $\lambda
_{\mathrm{S}}$ and that it performs better than the rest of the lower bounds.


\section{Discussion and conclusion\label{sec:discussion}}

Except in very particular cases, the SGR of matrix population models subjected
to environmental stochasticity can not be derived analytically. In this work
we have introduced four new upper and lower bounds for the SGR. These bounds,
which at least for two-stage models are easy derived as close form
expressions, can provide both quantitative and qualitative information about
ecological models. Specially relevant is the fact that they can be used to
give sufficient conditions for population growth or population extinction (see
the Introduction and Section \ref{sec:delta}).

The bounds for the SGR are been obtained following two different strategies.
The first one, valid for all kind of models, is a perturbation approach. The
second one, useful for Leslie models with two age classes, makes use of
certain bounds for the vector of population structure that we have introduced
and that generalize previous results in the field.

The different bounds have been analyzed through computer simulations in the
context of Leslie models with two age classes. In general, the best upper
bound does a better job at approximating $\lambda_{\mathrm{S}}$ than the best
lower bounds. As far as the upper bounds are concerned, the rate of growth
$\mu$ of the mean population size is the best upper bound for $\lambda
_{\mathrm{S}}$ in the vast majority of cases, and so it is the best quick
choice for the job. As far as the lower bounds are concerned, there is no
clear winner and the best bound depends strongly on the specific parameter
values, although a number of general guidelines have been provided. It has
also been found that Tuljapurkar's second order approximation to
$\lambda_{\mathrm{S}}$ underestimates $\lambda_{\mathrm{S}}$ in most (but not
all) cases and that it usually performs better that the lower bounds that have
been studied in this work.

Section \ref{sec:delta} illustrates one of the possible applications of the
bounds to the study of stochastic population models. There are many more
applications in which the bounds can be useful, and we proceed to give some
ideas about one of them. Let us consider an age-structured stochastic model
like that of Section \ref{sec:leslie} with two environments. In the first
environment the dominant eigenvalue of the corresponding PPM is greater that
one whereas in the case of the second environment it is less that one. By
using one lower bound and proceeding like in Section \ref{sec:delta}, we can
find a value $\pi^{\ast}\in(0,1)$ such that if the asymptotic probability
$\pi_{1}$ of environment one verifies $\pi_{1}>\pi^{\ast}$ then the population
grows to infinity. Similarly, the use of an upper bound allows one to find a
value $\hat{\pi}\in(0,1)$ such that the population goes extinct whenever
$\pi_{1}<\hat{\pi}$.

\section{}

\appendix

\section{Matrix models with environmental stochasticity\label{app1}}

The following theorem summarizes some results in the literature
\citep{furstenberg1960products, cohen1977ergodicity, tuljapurkar1980population}
about the dynamics of models of the kind (\ref{e1}) which are relevant in this work:

\begin{theorem}
\label{fust}Let us consider system (\ref{e1}) and assume that:

a. The set $\{\boldsymbol{A}_{\eta},\ \eta\in\mathcal{I}\}$\ of vital rates
matrices is ergodic, i.e., there exists a positive integer $g$ such that any
product of $g$ matrices (with repetitions allowed) drawn from $\mathcal{A}$ is
a positive matrix (i.e., its components are all positive) and, moreover, there
exist constants $\alpha,\beta>0$ such that for all $\eta\in\mathcal{I}$
\[
\max(\boldsymbol{A}_{\eta})\leq\beta\ \ ;\ \ \min\nolimits^{+}(\boldsymbol{A}%
_{\eta})\geq\alpha,
\]
where $\max(\boldsymbol{A}_{\eta})$ is the maximum of the elements of
$\boldsymbol{A}_{\eta}$ and $\min^{+}(\boldsymbol{A}_{\eta})$ is the minimum
of the positive elements of $\boldsymbol{A}_{\eta}$. Clearly, when
$\mathcal{I}$\ is finite, the second condition can be omitted.

b. $\tau_{t}$\ is a homogeneous Markov chain that has an unique stationary
distribution, and the transition probabilities of $\tau_{t}$ converge
geometrically to it. When $\mathcal{I}$ is finite, a sufficient condition for
this is that $\tau_{t}$ be an homogeneous irreducible and aperiodic
homogeneous Markov chain.

Then:

1. The joint process $\left(  \tau_{t},\boldsymbol{z}(t)/\left\Vert
\boldsymbol{z}(t)\right\Vert \right)  $ is a Markov chain that converges to a
stationary distribution $G$.

2. We can define the \textit{stochastic growth rate} (SGR) $\lambda
_{\mathrm{S}}$ of system (\ref{e1}) through $\log\lambda_{\mathrm{S}%
}:=\underset{t\rightarrow\infty}{\lim}\log\left\Vert \boldsymbol{z}%
(t)\right\Vert /t$ where the limit holds with probability one. Moreover,
$\lambda_{\mathrm{S}}$ is finite, is independent of the initial probabilities
of $\tau_{t}$ and of the initial (non-zero) population vector $\boldsymbol{z}%
_{0}\geq\mathbf{0}$ and can be calculated through
\begin{equation}
\log\lambda_{\mathrm{S}}=\mathbb{E}_{G}\log\left\Vert \boldsymbol{A}%
_{\tau_{t+1}}\frac{\boldsymbol{z}(t)}{\left\Vert \boldsymbol{z}(t)\right\Vert
}\right\Vert , \label{pepito}%
\end{equation}
where the expectation is taken with respect to the stationary distribution $G$.
\end{theorem}

\section{Proofs\label{sec:proofs}}

\noindent\textbf{Proof of the bounds (\ref{CotasPertBis})}. Following
\citep{furstenberg1960products, tuljapurkar1980population} the SGR of system
(\ref{e1}) can be calculated through any of the components of the matrix
product, i.e., for all $i,j=1,...,n$ we have%
\begin{equation}
\log\lambda_{\mathrm{S}}=\lim_{t\rightarrow\infty}\log\left(  \boldsymbol{A}%
_{\tau_{t}}\boldsymbol{A}_{\tau_{t-1}}\cdots\boldsymbol{A}_{\tau_{2}%
}\boldsymbol{A}_{\tau_{1}}\right)  ^{ij}/t, \label{Prod}%
\end{equation}
with probability one.

From (\ref{pert02}) we can write $A_{\eta}^{ij}=\left(  1+W_{\eta}%
^{ij}\right)  B_{\eta}^{ij}$ for all $\eta\in\mathcal{I}$ and $i,j=1,...,n,$
and therefore for all $i,j=1,...,n$ one has
\begin{gather*}
\left(  \boldsymbol{A}_{\tau_{t}}\boldsymbol{A}_{\tau_{t-1}}\cdots
\boldsymbol{A}_{\tau_{2}}\boldsymbol{A}_{\tau_{1}}\right)  ^{ij}\leq\\
\leq\left(  1+W_{\tau_{t}}^{M}\right)  \left(  1+W_{\tau_{t-1}}^{M}\right)
\cdots\left(  1+W_{\tau_{2}}^{M}\right)  \left(  1+W_{\tau_{1}}^{M}\right)
\left(  \boldsymbol{B}_{\tau_{t}}\boldsymbol{B}_{\tau_{t-1}}\cdots
\boldsymbol{B}_{\tau_{2}}\boldsymbol{B}_{\tau_{1}}\right)  ^{ij},
\end{gather*}
from where it follows%
\begin{gather}
\frac{\log\left(  \boldsymbol{A}_{\tau_{t}}\boldsymbol{A}_{\tau_{t-1}}%
\cdots\boldsymbol{A}_{\tau_{2}}\boldsymbol{A}_{\tau_{1}}\right)  ^{ij}}{t}%
\leq\label{pert06}\\
\frac{\log\left[  \left(  1+W_{\tau_{t}}^{M}\right)  \left(  1+W_{\tau_{t-1}%
}^{M}\right)  \cdots\left(  1+W_{\tau_{2}}^{M}\right)  \left(  1+W_{\tau_{1}%
}^{M}\right)  \right]  }{t}+\frac{\log\left(  \boldsymbol{B}_{\tau_{t}%
}\boldsymbol{B}_{\tau_{t-1}}\cdots\boldsymbol{B}_{\tau_{2}}\boldsymbol{B}%
_{\tau_{1}}\right)  ^{ij}}{t}.\nonumber
\end{gather}

Let us consider the scalar system%
\begin{equation}
x(t+1)=\left(  1+W_{\tau_{t+1}}^{M}\right)  x(t). \label{ESC}%
\end{equation}
Since, according to (\ref{pert04}), we have $1+W_{\eta}^{M}>0$ for all
$\eta\in\mathcal{I}$, the set $\left\{  1+W_{\eta}^{M},\eta\in\mathcal{I}%
\right\}  $ is ergodic and therefore system (\ref{ESC}) has a SGR. In order to
calculate it, we use (\ref{pepito}) and take into account that, since the
system is scalar, the process $x(t)/\left\Vert x(t)\right\Vert $ is trivial
(equals 1 with probability one). So the joint distribution of $\left(
\tau_{t},x(t)/\left\Vert x(t)\right\Vert \right)  $ coincides with the
distribution of $\tau_{t}$ and so
\begin{equation}
\lim_{t\rightarrow\infty}\frac{\log\left[  \left(  1+W_{\tau_{t}}^{M}\right)
\left(  1+W_{\tau_{t-1}}^{M}\right)  \cdots\left(  1+W_{\tau_{2}}^{M}\right)
\left(  1+W_{\tau_{1}}^{M}\right)  \right]  }{t}=\underset{\eta\in\mathcal{I}%
}{\sum}\pi_{\eta}\log\left(  1+W_{\eta}^{M}\right)  . \label{pert05}%
\end{equation}

From (\ref{Prod}) it follows that $\lim_{t\rightarrow\infty}\log\left(
\boldsymbol{A}_{\tau_{t}}\boldsymbol{A}_{\tau_{t-1}}\cdots\boldsymbol{A}%
_{\tau_{2}}\boldsymbol{A}_{\tau_{1}}\right)  ^{ij}/t=\log\lambda_{\mathrm{S}}$
and $\lim_{t\rightarrow\infty}\log\left(  \boldsymbol{B}_{\tau_{t}%
}\boldsymbol{B}_{\tau_{t-1}}\cdots\boldsymbol{B}_{\tau_{2}}\boldsymbol{B}%
_{\tau_{1}}\right)  ^{ij}/t=\log\lambda_{\mathrm{S}}^{\prime}$. Using this
facts, together with (\ref{pert05}), in (\ref{pert06}) we finally obtain%
\[
\log\lambda_{\mathrm{S}}\leq\log\lambda_{\mathrm{S}}^{\prime}+\underset
{\eta\in\mathcal{I}}{\sum}\pi_{\eta}\log\left(  1+W_{\eta}^{M}\right)  .
\]
The other inequality in (\ref{CotasPertBis}) follows from an analogous
reasoning substituting the upper bound in (\ref{pert06}) by the corresponding
lower bound.

\noindent\textbf{Proof of bounds (\ref{cotitas}) for population structure}. To
start with, since we are only interested in age structure, in matrices
$\boldsymbol{A}_{\eta}$ we can factor out the fertility of the young, i.e.,
\[
\boldsymbol{A}_{\tau_{t}}=f_{\tau_{t}}%
\begin{pmatrix}
1 & \varepsilon_{t}\\
\gamma_{t} & 0
\end{pmatrix}
,
\]
where we have defined%
\[
\varepsilon_{t}:=F_{\tau_{t}}/f_{\tau_{t}},\ \gamma_{t}:=s_{\tau_{t}}%
/f_{\tau_{t}},
\]
and we proceed to work with matrices
\[%
\begin{pmatrix}
1 & \varepsilon_{t}\\
\gamma_{t} & 0
\end{pmatrix}
.
\]

We use the same reasoning of Tuljapurkar in \citep{Tulja90}, and the
equivalent to his equation (8.1.3) reads
\[
\frac{z^{2}(t+1)}{z^{1}(t+1)}=\frac{\gamma_{t+1}}{1+\varepsilon_{t+1}%
\frac{z^{2}(t)}{z^{1}(t)}}=:f_{t+1}(z^{2}(t)/z^{1}(t)),
\]
where function $f_{t+1}$ is increasing as a function of $\gamma_{t+1}$,
decreasing as a function of $\varepsilon_{t+1}$ and decreasing as a function
of $z^{2}(t)/z^{1}(t)$. From the first two assertions if follows that, if we
define for each value of $z^{2}(t)/z^{1}(t)$
\[
f_{\min}(z^{2}(t)/z^{1}(t)):=\frac{\gamma_{m}}{1+\varepsilon_{M}\frac
{z^{2}(t)}{z^{1}(t)}},\text{ }f_{\max}(z^{2}(t)/z^{1}(t)):=\frac{\gamma_{M}%
}{1+\varepsilon_{m}\frac{z^{2}(t)}{z^{1}(t)}},
\]
one has
\[
f_{\min}(z^{2}(t)/z^{1}(t))\leq f_{t+1}(z^{2}(t)/z^{1}(t))\leq f_{\max}%
(z^{2}(t)/z^{1}(t)).
\]
From here, and reasoning as in \citep[Section 8.1]{Tulja90}, we obtain that,
for any initial vector and large enough $t$, the expression $\delta\leq
z^{2}(t)/z^{1}(t)\leq\kappa$ holds where
\begin{equation}
\delta=\lim_{t\rightarrow\infty}\delta(t),\ \ \kappa=\lim_{t\rightarrow\infty
}\kappa(t), \label{ppl}%
\end{equation}
and $\delta(t)$ and $\kappa(t)$ are defined iteratively through%
\begin{align}
\delta(1)  &  =(f_{\min}\circ f_{\max})(z^{2}(0)/z^{1}(0)),...,\ \delta
(t)=(f_{\min}\circ f_{\max})(\delta(t-1))\label{ppl1}\\
\kappa(1)  &  =(f_{\max}\circ f_{\min})(z^{2}(0)/z^{1}(0)),...,\ \kappa
(t)=(f_{\max}\circ f_{\min})(\kappa(t-1)).\nonumber
\end{align}

The dynamics of $\delta(t)$ can be studied as follows: Let
\[
\boldsymbol{X}_{\min}=%
\begin{pmatrix}
1 & \varepsilon_{M}\\
\gamma_{m} & 0
\end{pmatrix}
,\text{ \ \ \ \ }\boldsymbol{X}_{\max}=%
\begin{pmatrix}
1 & \varepsilon_{m}\\
\gamma_{M} & 0
\end{pmatrix}
,
\]
and take any positive initial vector $\boldsymbol{z}(0).$ If we apply matrix
$\boldsymbol{X}_{\max}$ to $\boldsymbol{z}(0)$ and then $\boldsymbol{X}_{\min
}$, the new population vector verifies that the ratio of adult to young is
given by $(f_{\min}\circ f_{\max})(z^{2}(0)/z^{1}(0))$. Therefore, from
(\ref{ppl}) and (\ref{ppl1}) we have $\delta=v^{2}/v^{1},$ where
$\boldsymbol{v}=(v^{1},v^{2})$ is the dominant eigenvector of matrix
\[
\boldsymbol{X}_{\min}\boldsymbol{X}_{\max}=%
\begin{pmatrix}
1+\varepsilon_{M}\gamma_{M} & \varepsilon_{m}\\
\gamma_{m} & \gamma_{m}\varepsilon_{m}%
\end{pmatrix}
.
\]
Straightforward computations show that $\delta$ has indeed the expression
given by (\ref{cotitas}).

Regarding $\kappa,$ the same reasonings above apply replacing product
$\boldsymbol{X}_{\min}\boldsymbol{X}_{\max}$ by $\boldsymbol{X}_{\max
}\boldsymbol{X}_{\min}$.

\noindent\textbf{Proof of bounds }(\ref{CotasCohenEstr}). We make use of
(\ref{pepito}) and therefore
\begin{align*}
\mathbb{E}_{G}\log\left\Vert \boldsymbol{A}_{\tau_{t+1}}\frac{\boldsymbol{z}%
(t)}{\left\Vert \boldsymbol{z}(t)\right\Vert }\right\Vert  &  =\mathbb{E}%
_{G}\log\left(  \overset{n}{\underset{i=1}{\sum}}\overset{n}{\underset
{j=1}{\sum}}A_{\tau_{t+1}}^{ij}\frac{z^{j}(t)}{\left\Vert \boldsymbol{z}%
(t)\right\Vert }\right)  \leq\mathbb{E}_{G}\log\left(  \overset{n}%
{\underset{i=1}{\sum}}\overset{n}{\underset{j=1}{\sum}}u_{\eta}^{j}%
A_{\tau_{t+1}}^{ij}\right)  =\\
&  =\mathbb{E}_{G}\log\left(  \overset{n}{\underset{j=1}{\sum}}u_{\eta}%
^{j}\overset{n}{\underset{i=1}{\sum}}A_{\tau_{t+1}}^{ij}\right)
=\mathbb{E}_{G}\log\left(  \overset{n}{\underset{j=1}{\sum}}u_{\eta}%
^{j}\varphi_{\tau_{t+1}}^{j}\right)  =\\
&  =\underset{\eta\in\mathcal{I}}{\sum}\pi_{\eta}\log\left(  \overset
{n}{\underset{j=1}{\sum}}u^{j}\varphi_{\eta}^{j}\right)  ,
\end{align*}
where in the last equality we have used that we are taking the expected value
of a random variable that is independent of the population structure
$\boldsymbol{z}(t)/\left\Vert \boldsymbol{z}(t)\right\Vert $. Thus, the
expression of $C_{IV}$ has been obtained. An analogous reasoning leads to the
expression of $c_{IV}$.


\smallskip

\smallskip

\smallskip

\smallskip

\smallskip

\noindent\textbf{Acknowledgements}

Author is supported by Ministerio de Econom\'{\i}a y Competitividad (Spain),
Project MTM2014-56022-C2-1-P.

\bibliographystyle{plain}
\bibliography{ConditionsGrowthAndExtinction}

\end{document}